\documentclass[amsmath,twocolumn,secnumarabic,%
    floatfix,amssymb,aps,nofootinbib,showkeys]{revtex4}
\usepackage{times,graphicx}
\usepackage{amsthm}
\usepackage[percent]{overpic}
\usepackage{placeins}
\usepackage{algorithm2e}
\usepackage{array}
\usepackage{booktabs}
\usepackage{epic}
\usepackage{eepic}
\usepackage[matrix,arrow]{xy}
\usepackage{url}

\newcommand{\tabnote}[1]{\caption{#1}}
\newcommand{\tbl}[2]{#2}

\renewcommand{\subsection}[1]{\medskip\noindent\textbf{#1.} }
\def\marginpar#1{}   

\let\lbl=\label
\def\label#1{\lbl{#1}\ifinner\else\marginpar{\ref{#1} #1}\ignorespaces\fi}

\newcommand{\cross}{\times}

\newcommand{\R}{\mathbb{R}}

\renewcommand{\phi}{\varphi}

\newcommand{\Thi}{\operatorname{Thi}}

\newcommand{\minRad}{\operatorname{minRad}}
\newcommand{\poca}{\operatorname{POCA}}

\newtheorem{theorem}{Theorem}[section]
\newtheorem{lemma}[theorem]{Lemma}
\newtheorem{definition}{Definition}

\graphicspath{{./figs/}}

{\makeatletter
 \gdef\xxxmark{%
   \expandafter\ifx\csname @mpargs\endcsname\relax 
     \expandafter\ifx\csname @captype\endcsname\relax 
       \marginpar{xxx}
     \else
       xxx 
     \fi
   \else
     xxx 
   \fi}
 \gdef\xxx{\@ifnextchar[\xxx@lab\xxx@nolab}
 \long\gdef\xxx@lab[#1]#2{{\bf [\xxxmark #2 ---{\sc #1}]}}
 \long\gdef\xxx@nolab#1{{\bf [\xxxmark #1]}}
}

\begin{document}

\title{A Fast Octree-Based Algorithm for Computing Ropelength}
\date{December 31, 2003; Revised: \today}

\author{Ted Ashton\footnote{email:ashted@uga.edu} and Jason
  Cantarella\footnote{email:cantarel@math.uga.edu}}
\affiliation{Department of Mathematics, University of Georgia, Athens, GA 30602}

\begin{abstract}
The \emph{ropelength} of a space curve is usually defined as the
quotient of its length by its \emph{thickness}: the diameter of the
largest embedded tube around the knot. This idea was extended to space
polygons by Eric Rawdon, who gave a definition of ropelength in terms
of doubly-critical self-distances (local minima or maxima of the distance
function on pairs of points on the polygon) and a function of the
turning angles of the polygon.
A naive algorithm for finding the doubly-critical self-distances of an
$n$-edge polygon involves comparing each pair of edges, and so takes
$O(n^2)$ time. In this paper, we describe an improved algorithm, based
on the notion of \emph{octrees}, which runs in $O(n \log n)$
time. The speed of the ropelength computation controls the performance
of ropelength-minimizing programs such as Rawdon and Piatek's TOROS.
An implementation of our algorithm is freely available under the GNU
Public License.
\end{abstract}


\maketitle

\section{Introduction}
For a $C^2$ curve in 3-space, \emph{ropelength} is the quotient of
the length of the curve by its \emph{thickness}: the diameter of the
largest embedded tube around the curve. Minimizing ropelength is the
same as fixing the diameter of the tube and minimizing its length---
if the tube is knotted, we are pulling the knot tight, and so the
minimum ropelength curves in any knot type are often called
\emph{tight} knots. Since the problem is such a natural one, the 
definition of thickness has been discovered and rediscovered by 
several authors\cite{MR42:8370,MR95k:58037,MR96d:58025}, with the
earliest results known (to these authors) on the problem credited 
to Kr{\"o}tenheerdt and Veit in 1976\cite{MR55:9069}.

In the past decade, there has been a great deal of interest in
exploring the geometry of tight knots; the definition of thickness has
been refined and fully understood\cite{MR2000c:57008}, it has been
shown that $C^{1,1}$ minimizers exist in each knot
type\cite{cks2,MR2003j:57010,MR2002m:74035}, some minimizing links have been
found\cite{cks2}, and a theory of ropelength criticality has started
to emerge\cite{cfksw,vdms}. The development of this theory has been
fueled by a steady stream of numerical data on ropelength minimizers,
from Pieranski's original SONO algorithm\cite{MR1702021} and Rawdon's
TOROS\cite{toros}, to second-generation efforts such as Smutny and
Maddocks' biarc computations\cite{clms,smbiarc} and the RIDGERUNNER
project of Cantarella, Piatek, and Rawdon. 
All of these algorithms have in their innermost loops a computation of
the ropelength of a curve in 3-space. 

Intuitively, the thickness of a
tube is controlled locally by the curvature of the core curve, and
globally by the approach of ``distant'' sections of the tube. 
Rawdon, in his thesis\cite{rawdonthesis}, defined a radius of curvature for a
corner of a polygon.  A given corner has two circles which are tangent to both
incident edges and tangent to one of the edges at its center.  He proved that
we can define a sensible polygonal radius of curvature as the radius of the
smaller of those two circles.  

More precisely:
\begin{definition}
If $P_n$ is a polygonal curve in $\R^3$ with edges $e_1, \dots, e_n$,
and $\alpha_i$ is the turning angle of the polygon made by edges
$e_i$ and $e_{i+1}$, then let
\begin{equation}
\minRad(P_n) = \min_{i \in 1, \dots, n} \left\{ \frac{|e_i|}{2
\tan\left( \frac{\alpha_i}{2} \right)}, \frac{|e_{i+1}|}{2 \tan\left(
\frac{\alpha_i}{2} \right)} \right\}
\end{equation}
where we take $e_{n+1} = e_{1}$ if the polygon is a closed curve, and
take $i \in 1, \dots, n-1$ otherwise.  
\end{definition}

\begin{definition}
Using the distance function on $P_n \cross P_n$ given by $D(x,y) = |x-y|$,
we say that a pair $xy$ of $P_n$ (bounding the chord $\overline{xy}$) is a
\emph{pair of closest approach} of $P_n$ if it is a non-trivial local minimum
of the distance function. The length of the shortest such chord is denoted
$\poca(P_n)$ (and we take $\poca(P_n) = \infty$ if no such chord exists).
\end{definition}

\begin{definition}
We define the \emph{thickness} of $P_n$ by
\begin{equation}
\Thi(P_n) = \min\left\{ 2\minRad(P_n), \poca(P_n) \right\}.
\end{equation}
\end{definition}
We note that the value which Rawdon uses in place of $\poca(P_n)$ in his
original definition of polygonal thickness\cite{rawdonthesis} is different.  In
particular, it is always finite.  But Rawdon reports that the equivalence of
the two definitions follows from results in an upcoming paper\cite{rawdon-pks}.

As computing the radius of curvature at a given corner only involves the edges
incident to that corner, computing $\minRad(P_n)$ requires only $O(n)$
time.  On the other hand, all previous efforts to compute thickness have used
some variant of Algorithm~\ref{alg:1} for computing $\poca(P_n)$.  This
algorithm is clearly $O(n^2)$.  So we have focused our attention on improving
the $\poca(P_n)$ calculation.

\begin{algorithm}
\SetLine
\SetInd{0.1in}{0.1in}
\For{$i = 1$ to $n$}{
	\For{$j = i+1$ to $n$}{
		check $e_i$ and $e_j$ for local min chords\;
		compare to previous shortest local min chord\;
	}
 }
 \caption{Standard Algorithm for $\poca(P_n)$.}
 \label{alg:1}
\end{algorithm}

Our algorithm concentrates on reducing the total number of edge-edge
checks performed by grouping the edges according to their positions in
space into a data structure known in computer graphics as an
\emph{octree}. We will use the octree to optimize the inner loop of
Algorithm~1, and show that we can isolate a constant-size set of
candidate $e_j$'s for any given $e_i$ in time $O(\log n)$. The new
algorithm will then perform $O(n \log n)$ edge-edge checks, and one
octree construction (which will also require time $O(n \log n)$). 

Before continuing, it is reasonable to ask whether such a complicated
algorithm can be implemented in a way that provides a practical
advantage over Algorithm~1. We believe that our implementation,
\texttt{liboctrope}, does. We give performance data for some test
problems in Section~\ref{sec:performance}. And more importantly, we invite
interested readers to download \texttt{liboctrope} and test the code
themselves (\url{http://ada.math.uga.edu/research/software/octrope}).

\section{Edge-Edge Checks}
\label{sec:ramps}

The quantity $\poca(P_n)$ is defined to be the smallest nontrivial local
minimum of the distance function $D(x,y)$ on pairs of points on the polygon
$P_n$. To understand it, we first make an observation about the nature
of these local minima.

\begin{lemma}
\label{lem:strut-types}
If we orient the curve $P_n$ and let $T^-(x)$, $T^+(x)$ denote the
inward and outward tangent vectors of $P_n$ at $x$ (they are different
if and only if $x$ is a vertex with nonzero turning angle). Every
pair $xy$ which locally minimizes $D:P_n \cross P_n \rightarrow \R$
has
\begin{eqnarray}
T^-(x) \cdot (y-x) \geq 0 & \quad & T^+(x) \cdot (y-x) \leq 0 \\
T^-(y) \cdot (x-y) \geq 0 & \quad & T^+(y) \cdot (x-y) \leq 0. 
\end{eqnarray}
We note that if $x$ is in the interior of an edge, then 
the above relations force $T^{\pm}(x)~\cdot~(x-y)~=~0$. 
\end{lemma}

\begin{figure}
\hspace*{\fill}
\begin{overpic}[scale=.95]{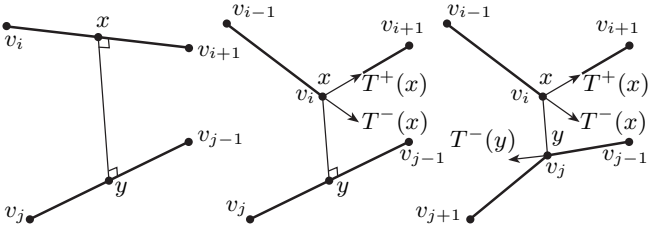}
\put(0.5,1){$v_j$}
\put(31,13){$v_{j-1}$}
\put(18,5){$y$}
\put(0.5,28){$v_i$}
\put(15,31){$x$}
\put(31,26.5){$v_{i+1}$}
\put(36.5,33){$v_{i-1}$}
\put(50,22){$x$}
\put(46.5,19.5){$v_i$}
\put(57,15){$T^-(x)$}
\put(57,21){$T^+(x)$}
\put(53,5){$y$}
\put(35.5,2){$v_j$}
\put(63,10){$v_{j-1}$}
\put(60,30.5){$v_{i+1}$}
\put(94,30.5){$v_{i+1}$}
\put(70.5,33){$v_{i-1}$}
\put(80.5,19.5){$v_i$}
\put(65.5,1){$v_{j+1}$}
\put(95,10){$v_{j-1}$}
\put(86,8){$v_j$}
\put(91.5,15){$T^-(x)$}
\put(92,21){$T^+(x)$}
\put(71,12){$T^-(y)$}
\put(85,22){$x$}
\put(87,13){$y$}
\end{overpic}
\hspace*{\fill}
\caption[Classifying Local Minima of $D$]
{We see the three cases in the proof of Lemma~\ref{lem:strut-types},
from left to right an \emph{edge-edge} pair, a \emph{vertex-edge} pair,
and a \emph{vertex-vertex} pair. In the center and right figures we see 
$T^{-}(x)$ and $T^+(x)$, and in the righthand figure we also see $T^{-}(y)$.
We have not drawn $T^+(y)$, but it would be colinear with $v_{j}v_{j+1}$
as with $T^+(x)$.}
\label{fig:strut-types}
\end{figure}

\begin{proof}
There are three cases: either both $x$ and $y$ are on the interior of
an edge, one is an edge point and one a vertex, or both are vertices,
as shown in Figure~\ref{fig:strut-types}. At $x$, the distance from
$y$ must not decrease to first order as one moves away from $x$ in
either direction along the curve: a computation verifies that this is
equivalent to the first line of the statement of the Lemma. A similar
argument at $y$ completes the proof.
\end{proof}

\begin{figure}
\hspace*{\fill}
\begin{overpic}{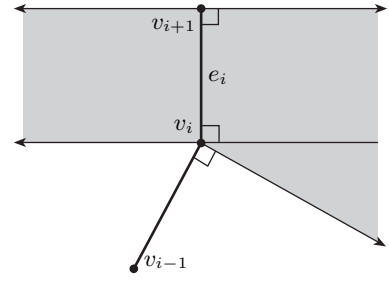}
\put(43,38){$v_{i}$} 
\put(37,65){$v_{i+1}$} 
\put(52,51){$e_i$}
\put(35,2){$v_{i-1}$}
\end{overpic}
\hspace*{\fill}
\caption[Possible locations for $y$ for $x$ on edge $e_i$.]  {The
shaded area represents the region of space in which the second point
$y$ of a locally minimal pair $xy$ can lie when $x$ is on the edge
$e_i$ or is the vertex $v_i$. This region consists of the infinite
slab of parallel planes normal to $e_i$ which pass through $e_i$,
together with the wedge extending from vertex $v_{i}$ in the outward
direction from the vertex.}
\label{fig:ramp}
\end{figure}

We now make a definition:
\begin{definition}
The $i^{\text{th}}$ \emph{ramp}, $R_i$ of a polygonal curve $P_n$ is the union
of the planes through edge $e_i = v_{i}v_{i+1}$ with normal vector $v_{i+1} -
v_i$, together with the wedge of vectors $w$ defined by the inequalities
\begin{equation}
(w - v_i) \cdot T^-(v_i) \geq 0 \quad (w - v_i) \cdot T^+(v_i) \leq 0.
\end{equation}
See Figure~\ref{fig:ramp}.
\end{definition}
This leads naturally to the Lemma:
\begin{lemma}
If $xy$ is a pair of points on $P_n$ which locally minimizes $D$, and
$x$ is on the half-open edge $e_i - \{v_{i+1}\}$, then $y$ is in the
ramp $R_i$.
\end{lemma}

\begin{proof}
If $x$ is in the interior of the edge, Lemma~\ref{lem:strut-types}
implies that $\overline{xy}$ must be perpendicular to $e_i$, and hence that $y$
is in the union of normal planes through $e_i$. If $x$ is at $v_i$,
the inequalities above are those of the statement of
Lemma~\ref{lem:strut-types}.
\end{proof}

\begin{figure}
\hspace*{\fill}
\begin{overpic}[scale=0.6]{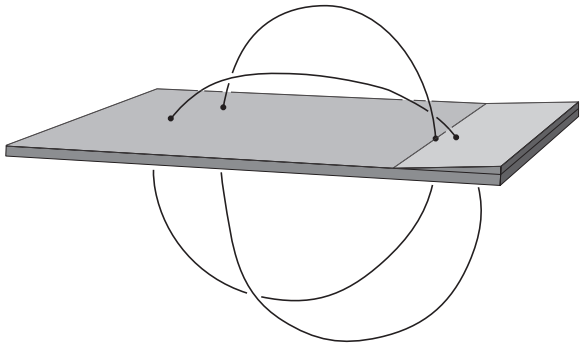}
\end{overpic}
\hspace*{\fill}
\caption[A typical ramp]{A typical ramp in a a trefoil of about 90 
edges consists of a very thin slab which only intersects the remainder
of the knot in a few places. If we could isolate ramp-knot intersections
quickly, it would reduce the number of edge-edge checks required to 
find the shortest $\poca$.}
\label{fig:thinslab}
\end{figure}

It may seem like we have only rephrased
Lemma~\ref{lem:strut-types}. In fact, we have gained an important
geometric insight about the problem-- for any edge $e_i$, the ramp
$R_i$ will probably be very close to a thin slab which only intersects
the remainder of $P_n$ in a few places (see
Figure~\ref{fig:thinslab}).  If we can isolate these intersections
quickly, we can complete the task of finding $\poca(P_n)$ by a more
detailed comparison of these candidates to $e_i$.

\section{The Octree Data Structure}
\label{sec:octree}

With the discussion in Section~\ref{sec:ramps}, we have reduced the
problem of identifying edges $e_j$ which may form locally minimal
pairs with points on edge $e_i$ to the problem of finding which edges
of $P_n$ intersect $e_i$'s ramp.  To do so efficiently, we will need a
new data structure for $P_n$: the octree\cite{jackinstanimoto}. 

\begin{figure*}[t]
\hspace*{\fill}
\begin{overpic}[scale=0.58]{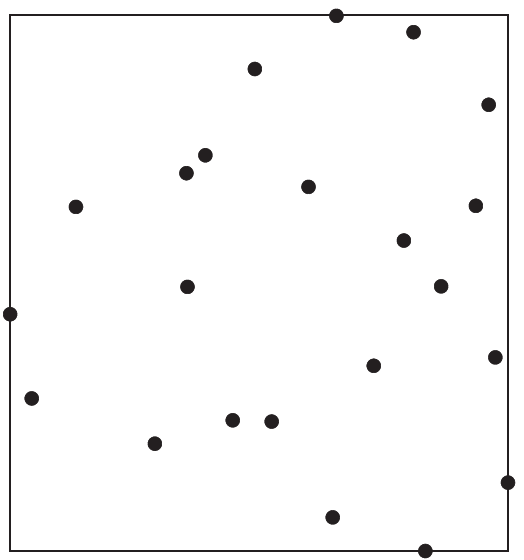}
\end{overpic}
\hspace*{\fill}
\begin{overpic}[scale=0.58]{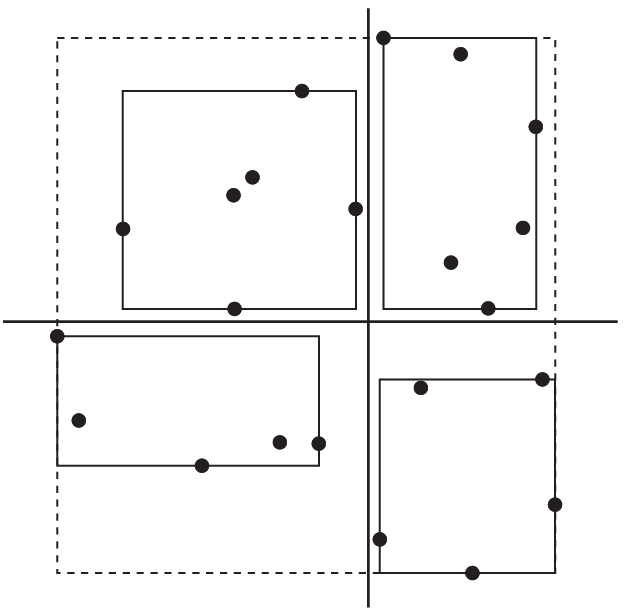}
\end{overpic}
\hspace*{\fill}
\begin{overpic}[scale=0.58]{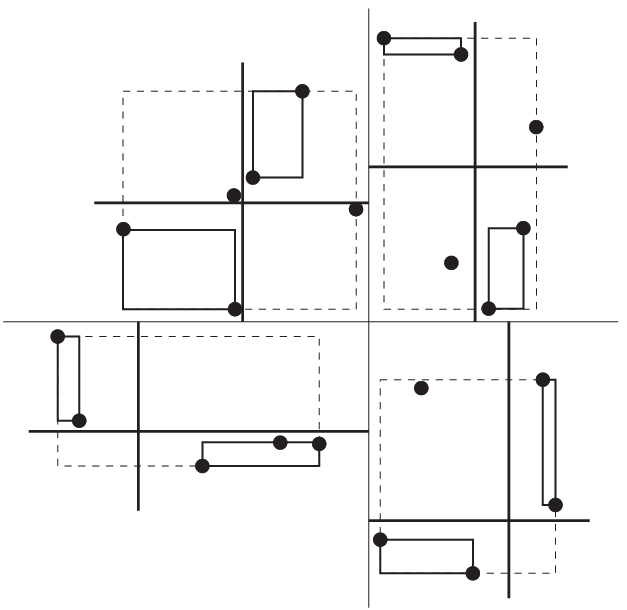}
\end{overpic}
\hspace*{\fill}
\caption[Construction of a Quadtree]{From left to right, these
pictures show three stages in the construction of a quadtree
representation (the planar version of an octree representation) of a
set of points. On the left, the bounding box of the entire point set
is computed. This is the root of the tree. In the center, we see the
points divided in two by $x$ and $y$ coordinates, and then grouped by
quadrant into four subcollections, with bounding boxes as shown. On
the right, we again divide the subcollections and group into
subquadrants. The resulting 3-level tree is shown in
Figure~\ref{fig:quadtree-rep}.}
\label{fig:quadtree-build}
\end{figure*}

\begin{figure*}[ht]

\begin{equation*}
\xymatrix@R=12pt@C=6pt{ 
 & & & & & & B \ar[dllll] \ar[dl] \ar[drr] \ar[drrrrr] \\ 
 & & B_1 \ar[dll] \ar[dl] \ar[d] \ar[dr] & & & B_2 \ar[dl] \ar[d] \ar[dr]
         \ar[drr] & & & B_3 \ar[d] \ar[dr] & & & B_4 \ar[dl] \ar[d] \ar[dr] \\ 
B_{11} \ar[d] & B_{12} \ar[d] & B_{13} \ar[d] & B_{14} \ar[d] & B_{21} \ar[d] &
B_{22} \ar[d] & B_{23} \ar[d] & B_{24} \ar[d] & B_{32} \ar[d] & B_{34} \ar[d] &
B_{41} \ar[d] & B_{42} \ar[d] & B_{43} \ar[d] \\ 
1 & 2 & 1 & 2 & 2 & 1 & 2 & 1 & 2 & 3 & 2 & 1 & 2 }
\end{equation*}

\caption[A Sample Quadtree]{This picture shows the quadtree
constructed in Figure~\ref{fig:quadtree-build} in a more familiar
form. The boxes are labelled according to the usual numbering
convention for quadrants of the plane, where the first quadrant is on
the top right and numbering proceedings counterclockwise.  The final
numbers show the number of points in each leaf box, and should be
compared to the boxes shown in Figure~\ref{fig:quadtree-build} in the
right-hand image.}
\label{fig:quadtree-rep}
\end{figure*}
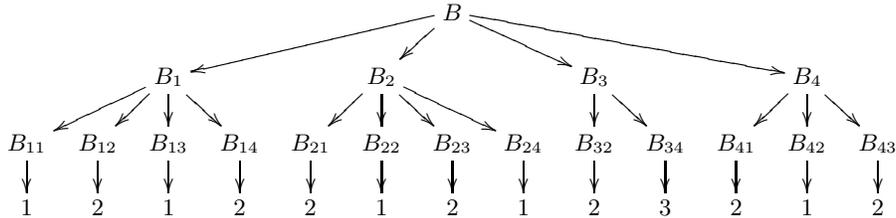

The octree representation of a collection of points in space is a tree
where each node represents the bounding box of a subset of that
collection. The eight daughter nodes of a parent represent the
bounding boxes of subsets of the points in the parent box created by
dividing that point set in two in each of the coordinate
directions. The most detailed octree representation of a point set has
leaf nodes which each contain a single point but it is common to
stop subdividing when the point sets are smaller than some fixed
number. Figure~\ref{fig:quadtree-build} illustrates three levels of
this process for a set of points in the plane, while
Figure~\ref{fig:quadtree-rep} shows the resulting tree.

From the description above, one can observe that it is easy to build an
octree using the recursive procedure of Algorithm~\ref{alg:buildoct}.

\begin{algorithm}
\SetLine
\SetInd{0.1in}{0.1in}
\AlgData{A node of the tree, a corresponding box $B$, the maximum number of
         points per box $m$ and a list of points.}
\AlgResult{An octree representation of this point list.}
\eIf{the list of points is no longer than $m$}{
	assign these points to this node\;
	make this node a leaf of the tree\;
	return\;
	}{
	partition $B$ into $8$ child boxes\;
	\For{each child box}{
		create a sublist of points intersecting that box\;
		recurse if this sublist is nonempty\;
		}
	}
\caption{One way to Build an Octree}
\label{alg:buildoct}
\end{algorithm}

For an $n$-point dataset, if one chooses each box partition so that no
child box contains more than half the total number of points in the
parent box, the number of levels in this tree is less than $\log_2 n$,
and one expects this algorithm to run in $O(n \log n)$ time. 
However, this algorithm involves a nontrivial amount of overhead in keeping
track of lists of points, and making procedure calls. Since we are
very concerned with the final performance of our implementation, we now
present a more insightful octree construction algorithm which has the 
same asymptotic time bound of $O(n \log n)$, but is much faster than
Algorithm~\ref{alg:buildoct} in practice.

To describe the new algorithm, we start with a numbering scheme. As 
we mentioned before, it is conventional to denote the upper right quadrant
of the plane by the number $1$, and proceed counterclockwise to the 
fourth quadrant on the lower right. A more natural numbering scheme 
assigns each quadrant a $2$-digit binary number, $d_x d_y$, where 
$d_x = 0$ has lower values of $x$ (the left hand side) and $d_x = 1$ has
higher values of $x$ (the right hand side), while $d_y=0$ denotes lower
values of $y$ (the bottom half), while $d_y=1$ denotes higher values
of $y$ (the top half). For octants in 3-space, we could assign three
digit binary numbers $d_x d_y d_z$ similarly.

Now consider the process of quadtree construction again. At the first
subdivision, we divide the point set in two parts by $x$-coordinate and 
by $y$-coordinate. This gives us four groups of points, which we can
number as above by the 2-digit binary numbers $d_x d_y$. These groups
are the members of the 4 boxes in the next level of the tree, as we saw 
above.

But there is something else to notice here: If the collection of
points is sorted by $x$ and by $y$, the digits $d_x$ and $d_y$ for any
particular point are the most significant binary digits of that
point's \emph{position} in the sorted array. Further, if we continue
to subdivide the points into fourths by $x$ and $y$, the next pair
of binary digits associated to each point, ${d_x}_1 {d_y}_1$ will be
the next pair of binary digits in that point's position in the $x$
and $y$ arrays as well. Again, for octrees the situation is similar,
but we sort by $z$ as well, and create a sequence of 3-digit binary
(or $1$-digit octal) numbers.

Continuing this process, we see that each point in the collection has
a unique \emph{octal tag} generated by interleaving the binary digits of its
position in the sorted $x$, $y$, and $z$ arrays. This tag specifies its
position in the octree.  Further, if we made a least-first traversal
of the octree (descending to octants in the order of their octal
labels), the order in which we would encounter the points would be by
increasing octal tags.  These observations give rise to a new
octree-building algorithm:

\begin{algorithm}
\SetLine
\SetInd{0.1in}{0.1in}
\AlgData{A list of points in $\R^3$.}
\AlgResult{An octree representation of point list.}
Sort the points by $x$, $y$, and $z$ coordinates\;
Shuffle binary digits of array positions to create octal tags\;
Sort again by octal tags\;

Build tree from this traversal-ordered list\;
\caption{A faster octree-building algorithm}
\label{alg:buildoct-fast}
\end{algorithm}

The problem of building a tree from a traversal-ordered list of its
contents is a standard one in computer science. Our particular
solution is discussed in some detail in
Section~\ref{sec:implementation} below. We note that building the
octree from the list also has time complexity $O(n \log n)$, since
every node in the octree must be visited, but that this algorithm is
still much faster than the previous method of octree construction
(Algorithm~\ref{alg:buildoct}).  We are among many rediscoverers of
this method of octree construction, which traces its roots to the
``linear quadtree'' construction of Gargantini\cite{gargantini}.

\section{The Core of the Algorithm}
\label{sec:core}

We can now describe our algorithm. Given any $e_i$, we must identify
all edges $e_j$ which might be part of a shortest $\poca$ with $e_i$. 
Such edges must obey two conditions: they must intersect $e_i$'s ramp,
and they must be closer to $e_i$ than the shortest $\poca$ found so 
far. Since both conditions can be checked for sub-boxes of the octree,
we can use them to eliminate groups of $e_j$ from consideration before
performing edge-edge checks.

In pseudo-code, this is a collection of $n$ calls to
the (recursive) Algorithm~\ref{alg:octrope} (one for each $e_i$).
We refer to the entire algorithm ($\minRad$ computation, octree 
construction by Algorithm~\ref{alg:buildoct-fast}, and calls to 
Algorithm~\ref{alg:octrope} for each edge) as {\bf Octrope}.

\begin{algorithm}
\SetLine
\SetInd{0.1in}{0.1in}

\AlgData{An octree node, the current minimum $\poca$ length $s$, the maximum
number of edges per box $m$ and a ramp from $e_i$.}
\AlgResult{All minimum-length $\poca$s between $e_i$ and edges in this subtree
  and (perhaps) a smaller value for the current minimum $\poca$ length $s$.}

\If{this box is within $s$ of $e_i$} {
\If{this box intersects the ramp from $e_i$}{
	\eIf{this box is a leaf}{
		check the (at most $m$) edges against $e_i$\;
		\If{$\poca$s of length $\leq s$ are found}{update $s$\; 
                  return list of minimum length $\poca$s\;}
		}{
		\For{each nonempty child node}{
			recurse on the child node\;
			}
		}
	
	}
}
\caption{Recursively identifying candidate $e_j$'s.}
\label{alg:octrope}
\end{algorithm}

Each call to this algorithm might require it to traverse the entire
depth of the octree before reaching leaf nodes and performing the
edge-edge checks.  Yet this depth is bounded above by $\log_2 n$, so
the expected running time for the algorithm is $O(\log n)$. In
pathological cases, many or all of the boxes may intersect the
ramp. If all the boxes intersect all the ramps, this algorithm may be
asymptotically slower than the naive one: we are forced to visit $O(n
\log n)$ tree nodes against each of $n$ edges, for a total time complexity of
$O(n^2 \log n)$. We have seen Algorithm~\ref{alg:1} outperform 
Algorithm~\ref{alg:octrope} only for a particularly bad class 
of examples: knots formed by connecting vertices chosen at
random inside a fixed volume. (See Section~\ref{sec:performance} for
details.)

\section{Implementation Issues}
\label{sec:implementation}

While being able to replace an $O(n^2)$ algorithm with one which is
$O(n \log n)$ will certainly save time \emph{for large enough values
of $n$}, there is no guarantee that this will help with problems of
practical size. Indeed, Algorithm~\ref{alg:octrope} threatens to
consume a fair amount of overhead, while Algorithm~\ref{alg:1}
involves only edge-edge checks, which could be coded very
efficiently. So in this section we turn our attention from the ``$n \log
n$'' to its multiplier --- from mathematics to program design.
In this discussion, we'll refer to function names and prototypes from
our publically available library version of {\bf Octrope}, which is
called \texttt{liboctrope}.

\subsection{The depth of the octree}
Since searching the octree involves some overhead, it is to be 
expected that we will not get the best performance from the deepest
octree.  Rather, we expect it to be more efficient to group some number of
edges in each box and do simple checks between the current edge and the
edges in an implicated leaf box. 

We implement this by bounding the maximum number of levels in the tree by some
$\ell$ and using that bound to calculate $m$, the maximum number of edges in
any leaf box, by the formula $m = \left\lceil \frac{n}{2^{\ell-1}}
\right\rceil$ (where $\lceil r \rceil$ is the least integer greater than or
equal to $r$).  The value of $\ell$ can be set by the user, using the
\texttt{octree\_set\_levels} call or it will default to $\ell = \left\lceil
\frac{3}{4} \log_2 n \right\rceil$, a formula at which we arrived empirically.

\subsection{A concrete example}
\begin{figure}
\hspace*{\fill}
\begin{overpic}{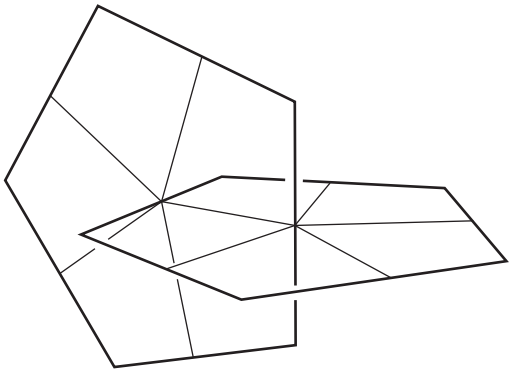}
\end{overpic}
\hspace*{\fill}
\caption[The pentagonal Hopf link]{Our example is a polygonal
  Hopf-link approximation composed of two regular pentagons. The
  lighter lines are the $9$ minimum length $\poca$s for which the {\bf
    Octrope} algorithm is searching.  They extend from the midpoints
  of one side of each interlocked polygon to all the midpoints of the
  edges of the other.}
\label{fig:pentahopf}
\end{figure}

We will now trace through our implementation of the fast octree
construction procedure of Algorithm~\ref{alg:buildoct-fast} for a
particular example: a Hopf link where each edge is given by a regular
pentagon (see Figure~\ref{fig:pentahopf} and Table~\ref{hopf}).  In this
example, $\ell$ would default to~$3$ and~$m$ would then also equal~$3$.

\begin{table}
\tbl{Hopf link vertices}{
\begin{tabular}{ll} \hline
\multicolumn{2}{c}{Component} \\
\multicolumn{1}{c}{1} & \multicolumn{1}{c}{2} \\ \hline
$v_{00} = (14.5,20,0)$    & $v_{10} = (0,0,14.5)$     \\
$v_{01} = (23.5,-7.6,0)$  & $v_{11} = (0,27.6,23.5)$  \\
$v_{02} = (0,-24.7,0)$    & $v_{12} = (0,44.7,0)$     \\
$v_{03} = (-23.5,-7.6,0)$ & $v_{13} = (0,27.6,-23.5)$ \\
$v_{04} = (-14.5,20,0)$   & $v_{14} = (0,0,-14.5)$    \\ \hline
\end{tabular}
\label{hopf}
}
\tabnote{The approximate vertices of the
pentagonal Hopf link shown in Figure~\ref{fig:pentahopf} are given
here. We have rounded the numbers to the nearest tenth to simplify the
table. This does not affect the octree-building procedure under
discussion, but would change the picture shown in the figure above.}
\end{table}

\newcommand{\byx}{\texttt{by\_x}}
\newcommand{\byy}{\texttt{by\_y}}
\newcommand{\byz}{\texttt{by\_z}}
\newcommand{\byoct}{\texttt{by\_oct}}

\subsection{Sorting edges}
The algorithm begins by gathering all of the edges into a single
$n$-element array which we call \byoct.  To avoid double-checking edge
pairs later on, we need each edge to ``belong'' to only one of the
leaf boxes in our tree.  So we identify each edge by its midpoint and
store that, as well as the edge's length and its starting vertex in
\byoct.

As we create \byoct, we also build \byx, \byy\ and \byz, three
$n$-element arrays of pointers to the elements of \byoct.  We then
sort these by $x$, $y$, and $z$ order. The result is shown in 
Table~\ref{sorted_edges}.

\begin{table*}
\tbl{The edges after sorting by $x$, $y$ and $z$}{
\begin{tabular}{lll} \hline
\byx & \byy & \byz \\\hline
$e_{03} : (-19,6.2,0)$ &
$e_{02} : (-11.75,-16.15,0)$ &
$e_{13} : (0,13.8,-19)$ \\
$e_{02} : (-11.75,-16.15,0)$ &
$e_{01} : (11.75,-16.15,0)$ &
$e_{12} : (0,36.15,-11.75)$ \\
$e_{04} : (0,20,0)$ &
$e_{14} : (0,0,0)$ & 
$e_{00} : (19,6.2,0)$ \\ \addlinespace[5pt]
$e_{10} : (0,13.8,19)$ &
$e_{03} : (-19,6.2,0)$ &
$e_{01} : (11.75,-16.15,0)$ \\
$e_{11} : (0,36.15,11.75)$ &
$e_{00} : (19,6.2,0)$ &
$e_{02} : (-11.75,-16.15,0)$ \\
$e_{12} : (0,36.15,-11.75)$ &
$e_{10} : (0,13.8,19)$ &
$e_{03} : (-19,6.2,0)$ \\ \addlinespace[5pt]
$e_{13} : (0,13.8,-19)$ &
$e_{13} : (0,13.8,-19)$ &
$e_{04} : (0,20,0)$ \\
$e_{14} : (0,0,0)$ & 
$e_{04} : (0,20,0)$ &
$e_{14} : (0,0,0)$ \\ 
$e_{01} : (11.75,-16.15,0)$ &
$e_{11} : (0,36.15,11.75)$ &
$e_{11} : (0,36.15,11.75)$ \\ \addlinespace[5pt]
$e_{00} : (19,6.2,0)$ &
$e_{12} : (0,36.15,-11.75)$ &
$e_{10} : (0,13.8,19)$ \\\hline
\end{tabular}
\label{sorted_edges}
}
\tabnote{In this table,
we see the $10$ edges of the pentagons in Table~\ref{hopf} sorted
by $x$, $y$, and $z$. The edges are sorted by their midpoints, and numbered
by the index of their first vertices. The spacing reminds us that since
$m = 3$, we are grouping the midpoints by threes when constructing boxes.}
\end{table*}

We divide \byx, \byz, and \byz\ into sections of $m$ points each
(shown in Table~\ref{sorted_edges} by spacing) and walk through them,
labeling the edges with the binary numbers of the sections in which
they lie in the following unusual fashion: if $x = x_1x_2\dotsb
x_{\ell-1}$, $y = y_1\dotsb y_{\ell-1}$ and $z = z_1\dotsb z_{\ell-1}$
are the respective box numbers and their binary representations, we
interleave those bits to produce a single octal number, $z_1 y_1 x_1
z_2 y_2 \dotsb y_{\ell-1} x_{\ell-1}$. This is the \emph{octal tag}
of Section~\ref{sec:octree} above
\footnote{To construct octal tags, we take a single pass
simultaneously through \byx, \byy, and \byz, starting with the second
box, which has binary tag $00\dotsb 01$.  As we walk the arrays, we
spread the bits of the box number apart (e.g.  $1101 \to 1001000001$)
using a lookup table similar to that of Shaffer\cite{shaffer}, shift
them left 1 bit for $y$ or 2 bits for $z$, and \texttt{OR} them with
the tag constructed so far.  The tags are thus built up over time and
guaranteed to be correct only when we reach the end of the pass.}.
We then sort \byoct\ by that octal tag, as shown in Table~\ref{octal_tags}.

\begin{table}
\tbl{The edges with their octal tags.}{
\begin{tabular}{cccccrr} \hline
edge & $x$-box & $y$-box & $z$-box & bits & octal & decimal\\\hline
$e_{02}$ & 0 & 0 & 1 & 000100 & $04_8$ &  4 \\
$e_{03}$ & 0 & 1 & 1 & 000110 & $06_8$ &  6 \\
$e_{00}$ & 3 & 1 & 0 & 001011 & $13_8$ & 11 \\
$e_{01}$ & 2 & 0 & 1 & 001100 & $14_8$ & 12 \\
$e_{12}$ & 1 & 3 & 0 & 010011 & $23_8$ & 19 \\
$e_{13}$ & 2 & 2 & 0 & 011000 & $30_8$ & 24 \\
$e_{10}$ & 1 & 1 & 3 & 100111 & $47_8$ & 39 \\
$e_{14}$ & 2 & 0 & 2 & 101000 & $50_8$ & 40 \\
$e_{04}$ & 0 & 2 & 2 & 110000 & $60_8$ & 48 \\
$e_{11}$ & 1 & 2 & 2 & 110001 & $61_8$ & 49 \\ \hline
\end{tabular}
}
\label{octal_tags}
\tabnote{This table contains the
edges with their box numbers in the $x$, $y$, and $z$ directions, the
binary numbers generated by interleaving the bits of these box
numbers, and the corresponding octal tags (in octal and decimal).  The data is
sorted by octal tag, and so appears in the same order in which it appears in
the \byoct\ array.  In principle, as many as $m$ edges can share an octal tag,
which means that they occupy the same leaf node of the resulting octree, but
this does not happen in our example.}
\end{table}

As we discussed in Section~\ref{sec:octree}, the sorted \byoct\
array is in the same order as that of a traversal of the
full octree.  

\subsection{Building the tree}
The actual building of the octree can be approached in various
directions. We could simply use the \byoct\ array with no futher indexing,
traversing it with binary searches (an approach which saves space at the
expense of time).  On the other hand, if we are to index it, we can build our
index in a top-down fashion, establishing the root node and building out to 
the leaf nodes.  We can build in a bottom-up fashion, partitioning off parts
of \byoct\ as leaf nodes and collecting them together in groups until we 
reach a single top node.  Or we can (and do) use a ``sideways'' or
``limb-by-limb'' approach.
We take an array of $\ell$ box pointers and on them build the ``left-hand
limb'' of the tree, all the way from the smallest numbered leaf box
down to the root box.  Each of the boxes knows its first edge in
\byoct\ and how many edges it has (which are grouped together thanks
to the octal sort).

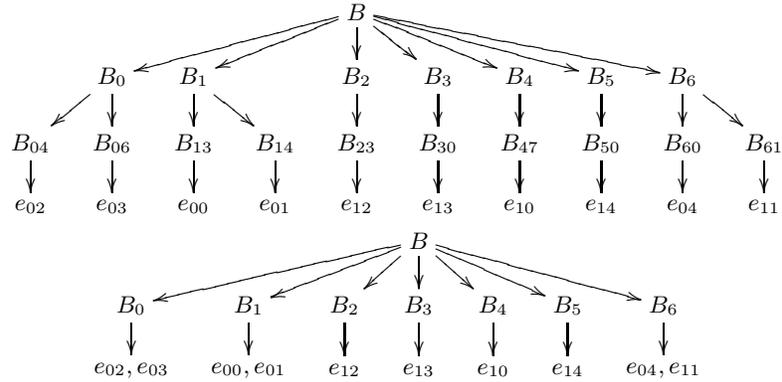
\begin{figure*}[b]
\begin{tabular}{c}
$$ \xymatrix@R=12pt@C=10pt{ & & & & B \ar[dlll] \ar[dll] \ar[d]
\ar[dr] \ar[drr] \ar[drrr] \ar[drrrr] & & & & & \\ & B_{0} \ar[dl]
\ar[d] & B_{1} \ar[d] \ar[dr] & & B_{2} \ar[d] & B_{3} \ar[d] & B_{4}
\ar[d] & B_{5} \ar[d] & B_{6} \ar[d] \ar[dr] & \\ B_{04} \ar[d] &
B_{06} \ar[d] & B_{13} \ar[d] & B_{14} \ar[d] & B_{23} \ar[d] & B_{30}
\ar[d] & B_{47} \ar[d] & B_{50} \ar[d] & B_{60} \ar[d] & B_{61} \ar[d]
\\ e_{02} & e_{03} & e_{00} & e_{01} & e_{12} & e_{13} & e_{10} &
e_{14} & e_{04} & e_{11} \\ }
$$\\
$$ \xymatrix@R=12pt@C=10pt{ & & & B \ar[dlll] \ar[dll] \ar[dl] \ar[d] \ar[dr]
\ar[drr] \ar[drrr] & & & \\ B_{0} \ar[d] & B_{1} \ar[d] & B_{2} \ar[d]
& B_{3} \ar[d] & B_{4} \ar[d] & B_{5} \ar[d] & B_{6} \ar[d] \\
e_{02},e_{03} & e_{00},e_{01} & e_{12} & e_{13} & e_{10} & e_{14} &
e_{04},e_{11} \\ }
$$
\end{tabular}
\caption[The octree without pruning (above) and with it (below)]{These two
trees show the octree as initially constructed (top) and after our pruning
procedure (bottom). This has grouped some edges together in single nodes (such
as $e_{02}$ and $e_{03})$ and deleted some nodes with only one child (such as
$B_{23}$). It is desirable to eliminate extra nodes of this kind, since even
though we keep the octal tree as compact in memory as possible, each jump from
node to node runs the risk of straying outside the memory cache of the
processor and incurring a delay as more information is loaded from main
memory.} 
\label{octree}
\end{figure*}

Then we walk once through \byoct, watching the octal tags.  As long as the
tag is the same as the one before it, we simply increment the count of
edges in that box.  When it changes, we do a binary \texttt{XOR} with
the previous tag to see how much they differ (that is, which of the
octal digits changed).  That tells how many of the boxes in this
``limb'' are complete.  After some cleanup (which may include pruning
the ``limb'') we leave those boxes and the create the new ones
necessary to hold this edge.  Figure \ref{octree} shows our example
tree after this process is complete.

\subsection{Searching the tree}
We have now created the tree and can move into using it.  We can check each
edge and its ramp against the tree, looking for leaf boxes on which to run
edge-edge checks.  Since $\poca$s are symmetric, we do not ever want to compare
the same pair of edges twice.  To avoid this we do the edge-edge check only if
the edge in question preceeds our chosen edge in \byoct.  By so doing, we can
eliminate entire boxes because their lowest edges are not in range.  The
improvement gained from technique has been significant.

\section{Performance}
\label{sec:performance}

We tested our algorithm using a $1.25$ Ghz $G4$ Macintosh computer on
high-resolution discretizations of trefoil knots and on (open) random walks.
We compiled our code with \texttt{gcc 3.3} and used the \texttt{-O3} option.
To make the tests, we compared the run times between \texttt{liboctrope} with
the tree depth set to $1$  and with the default tree depth of $\ell =
\left\lceil \frac{3}{4} \log_2 n \right\rceil$.  When the depth is 1, the
octree consists of a single box and $\frac{n(n-3)}{2}$ edge-edge checks are
performed during a run.  This turns \texttt{liboctrope} into a fairly 
efficient implementation of Algorithm~\ref{alg:1}.

\begin{figure*}[ht]

%
%
%
\setlength{\unitlength}{0.106pt}
\begin{picture}(3000,1800)(40,0)
\footnotesize

%
%

\put(40,820){\rotatebox{90}{Time (sec)}}
\put(1300,-40){Number of edges $n$}

\thicklines \path(287,166)(328,166)
\thicklines \path(2876,166)(2835,166)
\put(246,166){\makebox(0,0)[r]{ 0}}
\thicklines \path(287,321)(328,321)
\thicklines \path(2876,321)(2835,321)
\put(246,321){\makebox(0,0)[r]{ 20}}
\thicklines \path(287,476)(328,476)
\thicklines \path(2876,476)(2835,476)
\put(246,476){\makebox(0,0)[r]{ 40}}
\thicklines \path(287,632)(328,632)
\thicklines \path(2876,632)(2835,632)
\put(246,632){\makebox(0,0)[r]{ 60}}
\thicklines \path(287,787)(328,787)
\thicklines \path(2876,787)(2835,787)
\put(246,787){\makebox(0,0)[r]{ 80}}
\thicklines \path(287,942)(328,942)
\thicklines \path(2876,942)(2835,942)
\put(246,942){\makebox(0,0)[r]{ 100}}
\thicklines \path(287,1097)(328,1097)
\thicklines \path(2876,1097)(2835,1097)
\put(246,1097){\makebox(0,0)[r]{ 120}}
\thicklines \path(287,1252)(328,1252)
\thicklines \path(2876,1252)(2835,1252)
\put(246,1252){\makebox(0,0)[r]{ 140}}
\thicklines \path(287,1408)(328,1408)
\thicklines \path(2876,1408)(2835,1408)
\put(246,1408){\makebox(0,0)[r]{ 160}}
\thicklines \path(287,1563)(328,1563)
\thicklines \path(2876,1563)(2835,1563)
\put(246,1563){\makebox(0,0)[r]{ 180}}
\thicklines \path(287,1718)(328,1718)
\thicklines \path(2876,1718)(2835,1718)
\put(246,1718){\makebox(0,0)[r]{ 200}}
\thicklines \path(287,166)(287,207)
\thicklines \path(287,1718)(287,1677)
\put(287,83){\makebox(0,0){ 0}}
\thicklines \path(546,166)(546,207)
\thicklines \path(546,1718)(546,1677)
\put(546,83){\makebox(0,0){ 2000}}
\thicklines \path(805,166)(805,207)
\thicklines \path(805,1718)(805,1677)
\put(805,83){\makebox(0,0){ 4000}}
\thicklines \path(1064,166)(1064,207)
\thicklines \path(1064,1718)(1064,1677)
\put(1064,83){\makebox(0,0){ 6000}}
\thicklines \path(1323,166)(1323,207)
\thicklines \path(1323,1718)(1323,1677)
\put(1323,83){\makebox(0,0){ 8000}}
\thicklines \path(1582,166)(1582,207)
\thicklines \path(1582,1718)(1582,1677)
\put(1582,83){\makebox(0,0){ 10000}}
\thicklines \path(1840,166)(1840,207)
\thicklines \path(1840,1718)(1840,1677)
\put(1840,83){\makebox(0,0){ 12000}}
\thicklines \path(2099,166)(2099,207)
\thicklines \path(2099,1718)(2099,1677)
\put(2099,83){\makebox(0,0){ 14000}}
\thicklines \path(2358,166)(2358,207)
\thicklines \path(2358,1718)(2358,1677)
\put(2358,83){\makebox(0,0){ 16000}}
\thicklines \path(2617,166)(2617,207)
\thicklines \path(2617,1718)(2617,1677)
\put(2617,83){\makebox(0,0){ 18000}}
\thicklines \path(2876,166)(2876,207)
\thicklines \path(2876,1718)(2876,1677)
\put(2876,83){\makebox(0,0){ 20000}}
\thicklines \path(287,166)(2876,166)(2876,1718)(287,1718)(287,166)
\put(2548,1636){\makebox(0,0)[r]{Default tree depth}}
\thinlines \path(2589,1636)(2794,1636)
\thinlines \path(288,166)(288,166)(290,166)(291,166)(292,166)(293,166)(295,166)(296,166)(297,166)(299,166)(300,166)(301,166)(303,166)(304,166)(305,166)(306,166)(308,166)(309,166)(310,166)(312,166)(313,166)(314,166)(315,166)(317,166)(318,166)(319,166)(321,166)(322,166)(323,166)(325,166)(326,166)(327,166)(328,166)(330,166)(331,166)(332,166)(334,166)(335,166)(336,166)(337,166)(339,166)(340,166)(341,166)(343,166)(344,166)(345,166)(347,166)(348,166)(349,166)(350,166)(352,166)
\thinlines \path(352,166)(352,166)(365,166)(378,166)(391,166)(404,166)(416,166)(429,167)(442,167)(455,167)(468,167)(481,167)(494,167)(507,167)(520,167)(533,167)(546,167)(675,167)(805,168)(934,169)(1064,170)(1193,170)(1323,171)(1452,172)(1582,173)(1711,174)(1840,175)(1970,176)(2099,177)(2229,178)(2358,179)(2488,179)(2617,181)(2747,181)(2876,183)
\put(313,166){\makebox(0,2){$\blacktriangle$}}
\put(468,167){\makebox(0,2){$\blacktriangle$}}
\put(675,167){\makebox(0,2){$\blacktriangle$}}
\put(805,168){\makebox(0,2){$\blacktriangle$}}
\put(934,169){\makebox(0,2){$\blacktriangle$}}
\put(1064,170){\makebox(0,2){$\blacktriangle$}}
\put(1193,170){\makebox(0,2){$\blacktriangle$}}
\put(1323,171){\makebox(0,2){$\blacktriangle$}}
\put(1452,172){\makebox(0,2){$\blacktriangle$}}
\put(1582,173){\makebox(0,2){$\blacktriangle$}}
\put(1711,174){\makebox(0,2){$\blacktriangle$}}
\put(1840,175){\makebox(0,2){$\blacktriangle$}}
\put(1970,176){\makebox(0,2){$\blacktriangle$}}
\put(2099,177){\makebox(0,2){$\blacktriangle$}}
\put(2229,178){\makebox(0,2){$\blacktriangle$}}
\put(2358,179){\makebox(0,2){$\blacktriangle$}}
\put(2488,179){\makebox(0,2){$\blacktriangle$}}
\put(2617,181){\makebox(0,2){$\blacktriangle$}}
\put(2747,181){\makebox(0,2){$\blacktriangle$}}
\put(2876,183){\makebox(0,2){$\blacktriangle$}}
\put(2691,1636){\makebox(0,2){$\blacktriangle$}}
\put(2548,1553){\makebox(0,0)[r]{Tree depth $1$}}
\thicklines \path(2589,1553)(2794,1553)
\thicklines \path(288,166)(288,166)(290,166)(291,166)(292,166)(293,166)(295,166)(296,166)(297,166)(299,166)(300,166)(301,166)(303,166)(304,166)(305,166)(306,166)(308,166)(309,166)(310,166)(312,166)(313,166)(314,166)(315,166)(317,166)(318,166)(319,166)(321,166)(322,166)(323,166)(325,166)(326,166)(327,166)(328,166)(330,166)(331,166)(332,166)(334,166)(335,166)(336,166)(337,166)(339,166)(340,166)(341,166)(343,166)(344,166)(345,166)(347,166)(348,166)(349,166)(350,167)(352,167)
\thicklines \path(352,167)(352,167)(365,167)(378,167)(391,167)(404,168)(416,168)(429,169)(442,169)(455,170)(468,170)(481,171)(494,172)(507,172)(520,173)(533,174)(546,175)(675,186)(805,205)(934,237)(1064,280)(1193,329)(1323,382)(1452,441)(1582,507)(1711,580)(1840,660)(1970,746)(2099,840)(2229,944)(2358,1051)(2488,1167)(2617,1291)(2747,1426)(2876,1568)
\put(313,166){\circle{20}}
\put(468,170){\circle{20}}
\put(675,186){\circle{20}}
\put(805,205){\circle{20}}
\put(934,237){\circle{20}}
\put(1064,280){\circle{20}}
\put(1193,329){\circle{20}}
\put(1323,382){\circle{20}}
\put(1452,441){\circle{20}}
\put(1582,507){\circle{20}}
\put(1711,580){\circle{20}}
\put(1840,660){\circle{20}}
\put(1970,746){\circle{20}}
\put(2099,840){\circle{20}}
\put(2229,944){\circle{20}}
\put(2358,1051){\circle{20}}
\put(2488,1167){\circle{20}}
\put(2617,1291){\circle{20}}
\put(2747,1426){\circle{20}}
\put(2876,1568){\circle{20}}
\put(2691,1553){\circle{20}}
\thicklines \path(287,166)(2876,166)(2876,1718)(287,1718)(287,166)
\end{picture}

\caption[Comparative Performance of Octrope]{This plot shows the time
  (in seconds) required to find the ropelength for trefoils as a function
  of the number of edges $n$ in the polygonal knot.  The timings were
  computed on a $1.25$ Ghz Macintosh $G4$ computer, and represent
  averages over $10$ runs (for times above one second), or $100$ runs
  (for times below one second).  The data marked with
  \makebox[5pt]{\setlength{\unitlength}{0.120450pt}\begin{picture}(5,12)(0,0)\put(1,14){\makebox(0,2){$\blacktriangle$}}\end{picture}}
  comes from \texttt{liboctrope} with the default tree depth of $\ell
  = \left\lceil \frac{3}{4} \log_2 n \right\rceil$, while the data
  marked with
  \makebox[5pt]{\setlength{\unitlength}{0.120450pt}\begin{picture}(5,12)(0,0)\put(1,14){\circle{20}}\end{picture}}
  comes from \texttt{liboctrope} with the tree depth set to $1$ to
  force $n(n-3)/2$ edge-edge checks.  }

\label{fig:trefoildata}
\end{figure*}
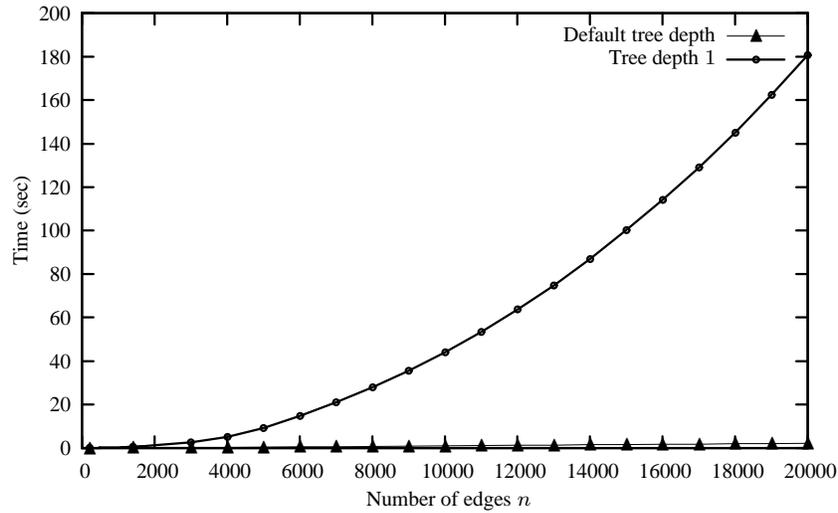

For both of these classes of knots, \textbf{Octrope} was much
faster than our reference implementation of Algorithm~\ref{alg:1}.
Figure~\ref{fig:trefoildata} shows the relative performance of the two
algorithms on trefoil knots given by $\gamma(\theta) = \left((1 +
\frac{2}{3} \cos 3\theta) \cos 2\theta, (1 + \frac{2}{3} \cos 3\theta)
\sin 2\theta,\frac{2}{3} \sin 3\theta\right)$. In fact, \textbf{Octrope}
outperformed Algorithm~\ref{alg:1} by an even greater margin on random walks.

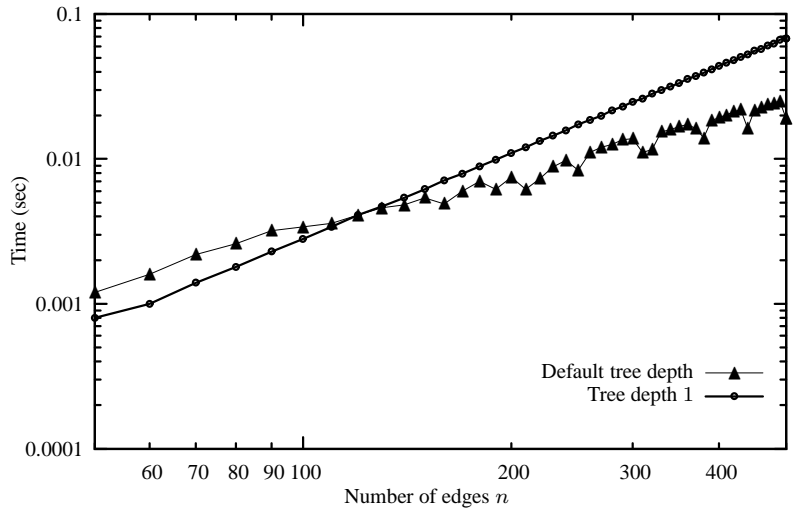
\begin{figure*}[t]


%
%

\setlength{\unitlength}{0.106pt}
\begin{picture}(3000,1800)(0,0)
\footnotesize


\put(110,820){\rotatebox{90}{Time (sec)}}
\put(1300,-30){Number of edges $n$}

\thicklines \path(410,166)(451,166)
\thicklines \path(2876,166)(2835,166)
\put(369,166){\makebox(0,0)[r]{ 0.0001}}
\thicklines \path(410,322)(430,322)
\thicklines \path(2876,322)(2856,322)
\thicklines \path(410,413)(430,413)
\thicklines \path(2876,413)(2856,413)
\thicklines \path(410,477)(430,477)
\thicklines \path(2876,477)(2856,477)
\thicklines \path(410,528)(430,528)
\thicklines \path(2876,528)(2856,528)
\thicklines \path(410,569)(430,569)
\thicklines \path(2876,569)(2856,569)
\thicklines \path(410,603)(430,603)
\thicklines \path(2876,603)(2856,603)
\thicklines \path(410,633)(430,633)
\thicklines \path(2876,633)(2856,633)
\thicklines \path(410,660)(430,660)
\thicklines \path(2876,660)(2856,660)
\thicklines \path(410,683)(451,683)
\thicklines \path(2876,683)(2835,683)
\put(369,683){\makebox(0,0)[r]{ 0.001}}
\thicklines \path(410,839)(430,839)
\thicklines \path(2876,839)(2856,839)
\thicklines \path(410,930)(430,930)
\thicklines \path(2876,930)(2856,930)
\thicklines \path(410,995)(430,995)
\thicklines \path(2876,995)(2856,995)
\thicklines \path(410,1045)(430,1045)
\thicklines \path(2876,1045)(2856,1045)
\thicklines \path(410,1086)(430,1086)
\thicklines \path(2876,1086)(2856,1086)
\thicklines \path(410,1121)(430,1121)
\thicklines \path(2876,1121)(2856,1121)
\thicklines \path(410,1151)(430,1151)
\thicklines \path(2876,1151)(2856,1151)
\thicklines \path(410,1177)(430,1177)
\thicklines \path(2876,1177)(2856,1177)
\thicklines \path(410,1201)(451,1201)
\thicklines \path(2876,1201)(2835,1201)
\put(369,1201){\makebox(0,0)[r]{ 0.01}}
\thicklines \path(410,1356)(430,1356)
\thicklines \path(2876,1356)(2856,1356)
\thicklines \path(410,1447)(430,1447)
\thicklines \path(2876,1447)(2856,1447)
\thicklines \path(410,1512)(430,1512)
\thicklines \path(2876,1512)(2856,1512)
\thicklines \path(410,1562)(430,1562)
\thicklines \path(2876,1562)(2856,1562)
\thicklines \path(410,1603)(430,1603)
\thicklines \path(2876,1603)(2856,1603)
\thicklines \path(410,1638)(430,1638)
\thicklines \path(2876,1638)(2856,1638)
\thicklines \path(410,1668)(430,1668)
\thicklines \path(2876,1668)(2856,1668)
\thicklines \path(410,1694)(430,1694)
\thicklines \path(2876,1694)(2856,1694)
\thicklines \path(410,1718)(451,1718)
\thicklines \path(2876,1718)(2835,1718)
\put(369,1718){\makebox(0,0)[r]{ 0.1}}
\thicklines \path(410,166)(410,186)
\thicklines \path(410,1718)(410,1698)
\thicklines \path(605,166)(605,186)
\thicklines \path(605,1718)(605,1698)
\put(605,83){\makebox(0,0){ 60}}
\thicklines \path(770,166)(770,186)
\thicklines \path(770,1718)(770,1698)
\put(770,83){\makebox(0,0){ 70}}
\thicklines \path(913,166)(913,186)
\thicklines \path(913,1718)(913,1698)
\put(913,83){\makebox(0,0){ 80}}
\thicklines \path(1040,166)(1040,186)
\thicklines \path(1040,1718)(1040,1698)
\put(1040,83){\makebox(0,0){ 90}}
\thicklines \path(1152,166)(1152,207)
\thicklines \path(1152,1718)(1152,1677)
\put(1152,83){\makebox(0,0){ 100}}
\thicklines \path(1895,166)(1895,186)
\put(1895,83){\makebox(0,0){ 200}}
\thicklines \path(1895,1718)(1895,1698)
\thicklines \path(2329,166)(2329,186)
\thicklines \path(2329,1718)(2329,1698)
\put(2329,83){\makebox(0,0){ 300}}
\thicklines \path(2637,166)(2637,186)
\thicklines \path(2637,1718)(2637,1698)
\put(2637,83){\makebox(0,0){ 400}}
\thicklines \path(2876,166)(2876,186)
\thicklines \path(2876,1718)(2876,1698)
\thicklines \path(410,166)(2876,166)(2876,1718)(410,1718)(410,166)
\put(2548,436){\makebox(0,0)[r]{Default tree depth}}
\thinlines \path(2589,436)(2794,436)
\thinlines \path(410,724)(410,724)(605,789)(770,860)(913,898)(1040,945)(1152,958)(1254,971)(1348,1000)(1433,1026)(1513,1036)(1587,1062)(1656,1040)(1721,1086)(1782,1121)(1840,1093)(1895,1136)(1947,1093)(1997,1130)(2044,1172)(2090,1196)(2134,1161)(2176,1224)(2216,1242)(2255,1253)(2293,1270)(2329,1273)(2364,1224)(2398,1234)(2431,1298)(2463,1305)(2494,1316)(2524,1324)(2554,1309)(2582,1276)(2610,1338)(2637,1348)(2663,1358)(2689,1369)(2714,1377)(2739,1310)(2763,1374)(2787,1384)(2810,1395)(2832,1401)(2854,1407)(2876,1345)(2876,1345)(2876,1345)
\put(410,724){\makebox(0,2){$\blacktriangle$}}
\put(605,789){\makebox(0,2){$\blacktriangle$}}
\put(770,860){\makebox(0,2){$\blacktriangle$}}
\put(913,898){\makebox(0,2){$\blacktriangle$}}
\put(1040,945){\makebox(0,2){$\blacktriangle$}}
\put(1152,958){\makebox(0,2){$\blacktriangle$}}
\put(1254,971){\makebox(0,2){$\blacktriangle$}}
\put(1348,1000){\makebox(0,2){$\blacktriangle$}}
\put(1433,1026){\makebox(0,2){$\blacktriangle$}}
\put(1513,1036){\makebox(0,2){$\blacktriangle$}}
\put(1587,1062){\makebox(0,2){$\blacktriangle$}}
\put(1656,1040){\makebox(0,2){$\blacktriangle$}}
\put(1721,1086){\makebox(0,2){$\blacktriangle$}}
\put(1782,1121){\makebox(0,2){$\blacktriangle$}}
\put(1840,1093){\makebox(0,2){$\blacktriangle$}}
\put(1895,1136){\makebox(0,2){$\blacktriangle$}}
\put(1947,1093){\makebox(0,2){$\blacktriangle$}}
\put(1997,1130){\makebox(0,2){$\blacktriangle$}}
\put(2044,1172){\makebox(0,2){$\blacktriangle$}}
\put(2090,1196){\makebox(0,2){$\blacktriangle$}}
\put(2134,1161){\makebox(0,2){$\blacktriangle$}}
\put(2176,1224){\makebox(0,2){$\blacktriangle$}}
\put(2216,1242){\makebox(0,2){$\blacktriangle$}}
\put(2255,1253){\makebox(0,2){$\blacktriangle$}}
\put(2293,1270){\makebox(0,2){$\blacktriangle$}}
\put(2329,1273){\makebox(0,2){$\blacktriangle$}}
\put(2364,1224){\makebox(0,2){$\blacktriangle$}}
\put(2398,1234){\makebox(0,2){$\blacktriangle$}}
\put(2431,1298){\makebox(0,2){$\blacktriangle$}}
\put(2463,1305){\makebox(0,2){$\blacktriangle$}}
\put(2494,1316){\makebox(0,2){$\blacktriangle$}}
\put(2524,1324){\makebox(0,2){$\blacktriangle$}}
\put(2554,1309){\makebox(0,2){$\blacktriangle$}}
\put(2582,1276){\makebox(0,2){$\blacktriangle$}}
\put(2610,1338){\makebox(0,2){$\blacktriangle$}}
\put(2637,1348){\makebox(0,2){$\blacktriangle$}}
\put(2663,1358){\makebox(0,2){$\blacktriangle$}}
\put(2689,1369){\makebox(0,2){$\blacktriangle$}}
\put(2714,1377){\makebox(0,2){$\blacktriangle$}}
\put(2739,1310){\makebox(0,2){$\blacktriangle$}}
\put(2763,1374){\makebox(0,2){$\blacktriangle$}}
\put(2787,1384){\makebox(0,2){$\blacktriangle$}}
\put(2810,1395){\makebox(0,2){$\blacktriangle$}}
\put(2832,1401){\makebox(0,2){$\blacktriangle$}}
\put(2854,1407){\makebox(0,2){$\blacktriangle$}}
\put(2876,1345){\makebox(0,2){$\blacktriangle$}}
\put(2876,1345){\makebox(0,2){$\blacktriangle$}}
\put(2691,436){\makebox(0,2){$\blacktriangle$}}
\put(2548,353){\makebox(0,0)[r]{Tree depth $1$}}
\thicklines \path(2589,353)(2794,353)
\thicklines \path(410,633)(410,633)(605,683)(770,759)(913,815)(1040,870)(1152,915)(1254,958)(1348,1000)(1433,1031)(1513,1062)(1587,1093)(1656,1124)(1721,1148)(1782,1174)(1840,1198)(1895,1222)(1947,1242)(1997,1265)(2044,1284)(2090,1302)(2134,1324)(2176,1340)(2216,1355)(2255,1374)(2293,1387)(2329,1405)(2364,1416)(2398,1434)(2431,1447)(2463,1459)(2494,1471)(2524,1486)(2554,1497)(2582,1509)(2610,1520)(2637,1533)(2663,1544)(2689,1553)(2714,1565)(2739,1574)(2763,1587)(2787,1593)(2810,1606)(2832,1613)(2854,1626)(2876,1631)(2876,1631)(2876,1631)
\put(410,633){\circle{20}}
\put(605,683){\circle{20}}
\put(770,759){\circle{20}}
\put(913,815){\circle{20}}
\put(1040,870){\circle{20}}
\put(1152,915){\circle{20}}
\put(1254,958){\circle{20}}
\put(1348,1000){\circle{20}}
\put(1433,1031){\circle{20}}
\put(1513,1062){\circle{20}}
\put(1587,1093){\circle{20}}
\put(1656,1124){\circle{20}}
\put(1721,1148){\circle{20}}
\put(1782,1174){\circle{20}}
\put(1840,1198){\circle{20}}
\put(1895,1222){\circle{20}}
\put(1947,1242){\circle{20}}
\put(1997,1265){\circle{20}}
\put(2044,1284){\circle{20}}
\put(2090,1302){\circle{20}}
\put(2134,1324){\circle{20}}
\put(2176,1340){\circle{20}}
\put(2216,1355){\circle{20}}
\put(2255,1374){\circle{20}}
\put(2293,1387){\circle{20}}
\put(2329,1405){\circle{20}}
\put(2364,1416){\circle{20}}
\put(2398,1434){\circle{20}}
\put(2431,1447){\circle{20}}
\put(2463,1459){\circle{20}}
\put(2494,1471){\circle{20}}
\put(2524,1486){\circle{20}}
\put(2554,1497){\circle{20}}
\put(2582,1509){\circle{20}}
\put(2610,1520){\circle{20}}
\put(2637,1533){\circle{20}}
\put(2663,1544){\circle{20}}
\put(2689,1553){\circle{20}}
\put(2714,1565){\circle{20}}
\put(2739,1574){\circle{20}}
\put(2763,1587){\circle{20}}
\put(2787,1593){\circle{20}}
\put(2810,1606){\circle{20}}
\put(2832,1613){\circle{20}}
\put(2854,1626){\circle{20}}
\put(2876,1631){\circle{20}}
\put(2876,1631){\circle{20}}
\put(2691,353){\circle{20}}
\thicklines \path(410,166)(2876,166)(2876,1718)(410,1718)(410,166)
\end{picture}

\caption[Crossover between Algorithm 1 and Octrope]{ This log-log plot
  shows the time (in seconds) required to find ropelength for trefoils
  as a function of the number of edges $n$ in the polygonal knot near
  the crossover point where \texttt{liboctrope} becomes faster than
  Algorithm~\ref{alg:1}. As in Figure~\ref{fig:trefoildata}, the
  timings were computed on a $1.25$ Ghz Macintosh $G4$ computer, and
  represent averages over $10$ runs (for times above one second), or
  $100$ runs (for times below one second).  The data marked with
  \makebox[5pt]{\setlength{\unitlength}{0.120450pt}\begin{picture}(5,12)(0,0)\put(1,14){\makebox(0,2){$\blacktriangle$}}\end{picture}}
  comes from \texttt{liboctrope} with the default tree depth of $\ell
  = \left\lceil \frac{3}{4} \log_2 n \right\rceil$, while the data
  marked with
   \makebox[5pt]{\setlength{\unitlength}{0.120450pt}\begin{picture}(5,12)(0,0)\put(1,14){\circle{20}}\end{picture}}
   comes from \texttt{liboctrope} with the tree depth set to $1$ to
   force $n(n-3)/2$ edge-edge checks.  The data shows that
  \texttt{liboctrope} is faster than our implementation of
  Algorithm~\ref{alg:1} for trefoil knots with more than about $120$
  edges.}

\label{fig:crossover}
\end{figure*}

For trefoils, Algorithm~\ref{alg:1} was sometimes faster than
\texttt{liboctrope} for very small numbers of edges. Figure~\ref{fig:crossover}
shows the performance of both algorithms near the crossover point in some
detail.

To understand the effect of varying the number of levels in the
octree, we also provide data for a 2499-edge random walk in
Table~\ref{tab:algorithms}.  Here we see a trade-off between edge-edge checks
and box/ramp checks as the octree resolution increases. Increasing the number
of levels in the octree from $9$ to $13$ cuts the number of final edge-edge
checks performed by a factor of $4$, but doubles the number of box/ramp checks.
Since the box/ramp checks are more computationally expensive, this is not a
favorable ratio, and the overall execution time increases. 

The data shows that we have been very effective at reducing the number of
edge-edge checks. On average, {\bf Octrope} compares each edge to less than
$13$ carefully chosen candidates when searching for minimum length $\poca$s.
\begin{table}
\tbl{Comparison of algorithms for random walk}{
\begin{tabular}{llclcl} \hline
                   & \multicolumn{5}{c}{Algorithm} \\ 
                   & \bf{Standard} & & \bf{Octrope}& & \bf{Max depth} \\ \hline
Octree levels      & $1$           & & $9 $        & &  $13$ \\
Edge-edge checks   & $3,121,251$   & & $32,033$    & &  $8189$ \\ 
Box/ramp checks    & $0$           & & $51,131$    & &  $93,187$ \\ 
Time               & $1.9$ sec    & &  $0.16$ sec & &  $0.22$ sec \\\hline
\end{tabular} 
\label{tab:algorithms}
}
\tabnote{This table compares the
performance of the \texttt{liboctrope} library on a $2500$-edge random walk
at three levels of tree depth: $1$ (the standard $O(n^2)$ algorithm), $\lceil
\frac{3}{4} \log_2 n \rceil$ ({\bf Octrope}), and $13$ (the maximum
resolution).}
\end{table}
 
It is worth noting that {\bf Octrope} is not guaranteed to outperform
the standard algorithm, even for large numbers of edges. For instance,
for random knots constructed by choosing vertices inside a fixed
volume, neither ramp-checking nor distance checking eliminates a
significant number of pairs from consideration. The edges are simply
too long, and pass too close to one another to decide in advance which
pairs are likely to control thickness. But even in this case,
\texttt{liboctrope} was only a few times slower than the standard
algorithm.
 
Currently, minimizing the ropelength of $850$-edge knots by simulated annealing
is a relatively taxing task, requiring a few weeks of computer time on a
standard desktop machine. Our timings above show that the ropelength
calculations involved in that process can be done 5 times faster using
\texttt{liboctrope} or that calculating the ropelength of a 2000-edge knot with
\texttt{liboctrope} would take the same time as finding that of the 850-edge
knot does now.

\section{Conclusions and Future Directions}
We have given an outline of an improved algorithm for computing the
ropelength of polygonal space curves in time $O(n \log n)$ and
contrasted it to the previous standard algorithm which required time
$O(n^2)$. We have implemented the algorithm efficiently in ANSI C, and 
given timings which show that our algorithm is also much faster in 
practice than previous methods used in the field.

The increase in speed from using our method should enable researchers
to consider significantly more complicated knots, and to get much
higher-resolution data for simpler knots. Both of these are valuable
goals. It has always been a goal of the geometric knot theory
community to apply our results to large biomolecules such as DNA and
proteins. Since these curves may involve thousands of vertices, they
have been out of the reach of tools based on Algorithm~\ref{alg:1}.

However, our methods do not entirely settle the problem of fast ropelength
computation. As mentioned in the introduction, Cantarella, Fu, Kusner,
Sullivan, and Wrinkle\cite{cfksw} have discovered tiny straight segments in a
ropelength-critical simple clasp.  These segments are a few one-thousands of
one unit in length out of a total clasp length of about $6$ units (a similar
clasp has been constructed by Starostin\cite{starostin}). To resolve these very
small scale phenomena numerically will require ropelength-minimized
configurations with tens of thousands of edges.

At about $3$ seconds per ropelength computation on a standard desktop machine,
it would simply be untenable to minimize the ropelength of a 20,000-edge knot
using our library and simulated annealing on a desktop machine. However,
{\bf Octrope} parallelizes well, so one could bring supercomputing cluster
machines to bear on the problem, reducing the time to evaluate a configuration
to tenths or hundredths of a second.  This might allow for a long enough
cooling schedule to resolve some small-scale phenomena, but there is no
guarantee.

We are hence considering two further approaches to the problem: the
use of Edelsbrunner's ``segment trees''\cite{edelsbrunner} and an
approach we call the {\bf Multiresolution Ropelength Algorithm}.
This algorithm is based on the idea that very high resolution
knots can be well-approximated by subsampling the vertex set. If one
keeps track of the distance between the subsampled knot and the
original, one can again eliminate groups of edges from edge-edge
checking. The potential advantages of this scheme are twofold: first,
the construction of the corresponding tree is linear in time and
second, the subsampled knots can be handled by {\bf Octrope}
itself. The disadvantage of the multiresolution algorithm is that it
will not help with random walks or other very complicated knots, such
as large protein backbones, as their subsamples will not be close to
the original curve.

We would also like to observe that while the discussion above is phrased in
terms of polygons, the general octree/ramp method is equally applicable to
other discretization schemes for curves, such as biarcs. In that case, the
relative speed advantage of this algorithm should be greater, since the
``edge-edge check'' for a pair of arcs or spline segments is much slower than
the edge-edge check for polygonal edges described above.

In conclusion, we hope that our algorithm and implementation will become a
standard software component in numerical investigations of the ropelength
problem. If others can improve our code, we hope that they will do so, and
invite them to contact us. We also hope that our public release of the library
(the first that we know of in geometric knot theory since Brakke's
Evolver\cite{MR93k:53006}) will inspire others in the field to contribute from
their personal and laboratory collections of code to the public domain.
Those interested in obtaining \texttt{liboctrope} can turn to
\url{http://ada.math.uga.edu/research/software/octrope/} for further
information.

\section{Acknowledgements}

The authors are grateful to many colleagues, including Herbert Edelsbrunner,
Mark Peletier, and John Sullivan, for discussions about ropelength and
algorithms.  The 2002-2003 VIGRE group in Geometric Knot Theory (in particular
Xander Faber, Chad Mullikin, and Nancy Wrinkle) contributed to our
understanding of the computational issues surrounding ropelength, and Monica
Shaw and Allison Diana, members of the 2003 Summer Undergraduate Research
Experience, worked on a prototype implementation of the octrope algorithm.
Michael Piatek and Eric Rawdon served as the \texttt{liboctrope} beta-test
team, as well as contributing insight on efficient code and library design.
The authors would also like to acknowledge the support of the National Science
Foundation through the University of Georgia VIGRE grant (DMS-00-89927),
DMS-99-02397 (to Cantarella), and DMS-02-04826 (to Cantarella and Fu).

\bibliography{drl,thick,octrope-paper} 
\vspace{-3ex}
\end{document}